\documentstyle[12pt,epsf,amssymb]{article}

\font\tenmsb=msbm10
\font\sevenmsb=msbm7
\font\fivemsb=msbm5
 
\newfam\msbfam 
\textfont\msbfam=\tenmsb
\scriptfont\msbfam=\sevenmsb
\scriptscriptfont\msbfam=\fivemsb
 
\def\Bbb#1{\fam\msbfam\relax#1}

\newtheorem{thm}{Theorem}[section]
\newtheorem{prop}[thm]{Proposition}
\newtheorem{cor}[thm]{Corollary}
\newtheorem{lem}[thm]{Lemma}
\newtheorem{conj}[thm]{Conjecture}
\newtheorem{exa}[thm]{Example}
\newtheorem{defn}[thm]{Definition}

\newtheorem{rem}[thm]{Remark}
\newtheorem{note}[thm]{Notation}
\newtheorem{alg}[thm]{Algorithm}

\newcommand{\ben}{\begin{enumerate}}
\newcommand{\een}{\end{enumerate}}
\newcommand{\ble}{\begin{lem}}
\newcommand{\ele}{\end{lem}}
\newcommand{\bth}{\begin{thm}}
\renewcommand{\eth}{\end{thm}}
\newcommand{\bpr}{\begin{prop}}
\newcommand{\epr}{\end{prop}}
\newcommand{\bco}{\begin{cor}}
\newcommand{\eco}{\end{cor}}
\newcommand{\bcon}{\begin{conj}}
\newcommand{\econ}{\end{conj}}
\newcommand{\bde}{\begin{defn}}
\newcommand{\ede}{\end{defn}}
\newcommand{\bex}{\begin{exa}}
\newcommand{\eex}{\end{exa}}
\newcommand{\brem}{\begin{rem}}
\newcommand{\erem}{\end{rem}}
\newcommand{\bnot}{\begin{note}}
\newcommand{\enot}{\end{note}}
\newcommand{\balg}{\begin{alg}}
\newcommand{\ealg}{\end{alg}}

\newcommand{\bib}{thebibliography}
\newcommand{\qed}{\square}

\newcommand{\C}{{\Bbb C}}

\newcommand{\PP}{{\Bbb P}}
\newcommand{\R}{{\Bbb R}}

\begin{document}

\title{Solving the Braid Word Problem Via the Fundamental Group} 

\author{S. Kaplan, M. Teicher$^1$ \\
\small Department of Mathematics and Computer Sciences\\
\small Bar-Ilan University\\
\small Ramat-Gan, Israel\\
\small \{kaplansh,teicher\}@macs.biu.ac.il}

\stepcounter{footnote}\footnotetext{Partially supported by the Emmy Noether Research Institute for
Mathematics, Bar-Ilan University and the Minerva Foundation, Germany.\\
Partially supported by the Excellency Center ``Group theoretic methods in
the study of algebraic varieties'' of the National Science Foundation of
Israel.}

\date{\today \\[1in]
}

\maketitle

\begin{abstract}
The word problem of a group is a very important question. The word problem in the braid 
group is of particular interest for topologists, algebraists and geometers. In \cite{WORD1}
we have looked at the braid group from a topological point of view, and thus using a new
computerized representation of some elements of the fundamental group we gave a solution
for its word problem. 
In this paper we will give an algorithm that will make it possible to transform the new
presentation from \cite{WORD1} into a syntactic presentation. This will make it possible
to computerize the group operation to sets of elements of the fundamental group, called a g-base, which are isomorphic
to the braid group. More over we will show that it is sufficient enough to look at the
syntactic presentation in order to solve the braid word problem, resulting with a better
and faster braid word solution.
\end{abstract}

\tableofcontents

\section*{Introduction}
Let $D$ be a closed disk, and $K=\{k_1,...,k_n\}$ be $n$ points in $D$. Let $B$ be 
the group of all diffeomorphisms $\beta$ of $D$ such that $\beta(K)=K$, and
$\beta |_{\partial D}={\rm Id} |_{\partial D}$. The braid group is derived from $B$ by identifying two 
elements if their actions on $\pi _1(D \setminus K,u)$ are equal. To simplify the algorithm, 
we choose a geometric base of $\pi _1(D \setminus K,u)$, and
we look at the action of $\beta \in B$ on the elements of this geometrical base. This yields a group
which is isomorphic to the braid group.

In \cite{WORD1} we gave a geometric approach for solving the braid word problem based on
this definition of the braid group. Now, we will show a method that will enable us to add
to the new computerized presentation of the braid group that was represented in
\cite{WORD1} a group structure, which means describing a way for computing the result of
multiplication of two braids only when given their two representations of the g-base.

This yields us with a new syntactic approach for encoding the elements of the g-base of
the fundamental group, that will make it possible to improve the running time of the
algorithm that solves the braid word problem, resulting with an even more practical
algorithm for short braid words.

In section 1, we will give the topological definition of the braid group, based on some 
facts concerning the fundamental group of a punctured disk, and in the end of this section we
will define the braid word problem. In section 2, we will recall the new presentation of
the elements of the g-base of the fundamental group that was presented in \cite{WORD1}.
In section 3 and 4, we will introduce the syntactic presentation of elements of the g-base and
the methods to transform the geometrical presentation into the syntactic and vice versa.
Section 5 will be dedicated to adding a group structure on the new presentations. Section
6 will show the improvement of the solution of the braid word problem that was given in
\cite{WORD1}, and section 7 will give conclusions, future applications of the
new presentation, and further plans.

\section{The topological definition of the braid group}
In the first part of this section we will give the definition of the fundamental group and some aspects of it,
while on the second part we will give the topological definition of the braid group which
is different then the widely used algebraic definition introduced by Artin \cite{Artin}.

\subsection{The fundamental group}

Let $D$ be a topological space. We fix a base point $u \in D$. A \underline{path} in $D$ is a
continuous map $\gamma:[0,1] \to D$. A \underline{loop} based at $u$ is a path on $D$ such that
$\gamma (0)=\gamma (1)=u$.
Two loops $\gamma _1$ and $\gamma _2$ are said to be \underline{homotopic} if there is a continuous
map $G:[0,1] \times [0,1] \to D$ such that $G(0,t)=\gamma _1(t)$ and $G(1,t)=\gamma _2(t)$
for all $0 \leq t \leq 1$, and $G(s,0)=G(s,1)=u$ for all $0 \leq s \leq 1$. 
Homotopy is an equivalence relation on the set of all loops based at $u$. 
\bde
\underline{The fundamental group} of $D$ is the set of homotopy classes of loops based at $u$, and is
denoted by $\pi _1(D,u)$. The operation of concatenation of loops forms a group structure
on $\pi _1(D,u)$.
\ede

We will reduce our look at the fundamental group to the case where $D$ is a closed disk,
and $K=\{k_1,.,,,k_n\}$ is a finite set such that $K \subset int(D)$. 
When we will refer to $D$ and $K$ from practical algorithmic point of view we will reduce 
our look even more to the case were $D$ is the closed unit disk and $K$ is a finite set of points 
ordered from left to right on the x-axis.

\brem
It is known that the fundamental group of a punctured disk with $n$ holes is a free group on $n$ generators.
\erem

We fix an orientation on $D$. Let $q$ be a simple path connecting $u$ with one of the $k_i$, say $k_{i_0}$, such that
$q$ does not meet any other point $k_j$ where $j \neq i_0$. To $q$ we will assign a
loop $l(q)$ (which is an element of $\pi _1(D \setminus K,u)$) as follows:

\bde{\underline{$l(q)$}}

Let $c$ be a simple loop equal to the (oriented) boundary of a small neighborhood $V$ of
$k_{i_0}$ chosen such that $q'=q \setminus (V \cap q)$ is a simple path.
Then \underline{$l(q)=q' \cup c \cup q'^{-1}$}. We will use the same notation for the element of $\pi
_1(D \setminus K,u)$ corresponding to $l(q)$.
\ede

\bde
Let  $(T_1,...,T_n)$ be an ordered set of simple paths in $D$ which connect the $k_i$'s with
$u$ such that:
\ben
\item
$T_i \cap k_j=\emptyset$ if $i \neq j$ for all $i,j=1,...,n$.
\item
$\displaystyle \bigcap_{i=1}^nT_i=\{u\}$.
\item
for a small circle $c(u)$ around $u$, each $u'_i=T_i \cap c(u)$ is a single point and the
order in $(u'_1,...,u'_n)$ is consistent with the positive orientation of $c(u)$.
\een

We say that two such sets $(T_1,...,T_n)$ and $(T'_1,...,T'_n)$ are \underline{equivalent} if 
$l(T_i)=l(T'_i)$ in $\pi _1(D \setminus K,u)$  for all $i=1,...,n$.

An equivalence class of such sets is called a \underline{bush} in $D \setminus K$.
\ede

\bde
A \underline{g-base} (geometrical base) of $\pi _1(D \setminus K,u)$ is an ordered 
free base of $\pi _1(D \setminus K,u)$
which has the form $(l(T_1),...,l(T_n))$, where $(T_1,...,T_n)$ is a bush in $D \setminus
K$.
\ede

We would like to point out a particular g-base which will be used in the paper.
Choose $T_i$ to be the straight line connecting $u$ with $k_i$, then we call $(l(T_1),...,l(T_n))$ the 
\underline{standard g-base of $\pi _1(D \setminus K,u)$} and it is shown in the following figure: 

\begin{center}
\epsfysize=4cm
\epsfbox{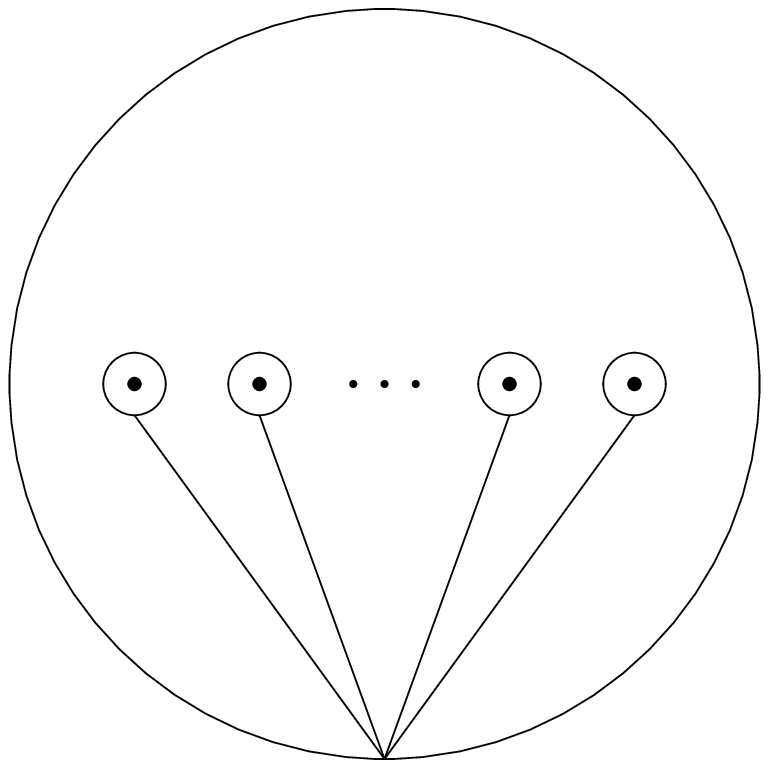}
\end{center}

We want to point out that the $n$ elements of the g-base generate the fundamental group of
the punctured disk $D \setminus K$, and we will call the elements of the standard g-base
the \underline{standard generators} of the fundamental group.

\subsection{The braid group}
Let $D,K,u$ be as above.

\bde
Let $B$ be the group of all diffeomorphisms $\beta$ of $D$ such that $\beta(K)=K$, $\beta
|_{\partial D}={\rm Id}|_{\partial D}$. For $\beta _1,\beta _2 \in B$ we say that $\beta _1$ is
\underline{equivalent} to $\beta _2$ if $\beta _1$ and $\beta _2$ induce the same automorphism of $\pi
_1 (D \setminus K,u)$. The quotient of $B$ by this equivalence relation is called 
\underline{the braid group $B_n[D,K]$} ($n=\#K$). The elements of $B_n[D,K]$ are called \underline{braids}.
\ede

\brem
For the canonical homomorphism $\psi:B \to Aut(\pi _1(D \setminus K,u))$
we actually have $B_n[D,K] \cong Im(\psi)$. 
\erem

We recall two facts from \cite{BGTI}[section III].

\ben
\item
If $K' \subset D'$, where $D$ is another disk, and $\#K=\#K'$ then $B_n[D,K] \cong B_n[D',K']$.
\item
Any braid $\beta \in B_n[D,K]$ transforms a g-base to a g-base. Moreover, 
for every two g-bases, there exists a unique braid which transforms one g-base to another.
\een

We distinguish some elements in $B_n[D,K]$ called \underline{half-twists}.

Let $D,K,u$ be as above. Let $a,b \in K$ be two points. We denote $K_{a,b}=K \setminus
\{a,b\}$. Let $\sigma$ be a simple path in $D \setminus (\partial D \cup K_{a,b})$ connecting
$a$ with $b$. Choose a small regular neighborhood $U$ of $\sigma$ and an orientation
preserving diffeomorphism $f:\R ^2 \to \C$ such that $f(\sigma )=[-1,1]$, $f(U)=\{z \in \C \ | \ |z| <2\}$. 

Let $\alpha (x)$, $0 \leq x$ be a real smooth monotone function such that:

$$\alpha (x)=\left \{ \matrix{1 & 0 \leq x \leq \frac{3}{2} \cr
                             0 & 2 \leq x} \right.$$

Define a diffeomorphism $h:\C \to \C$ as follows: for $z=re^{i\varphi} \in \C$ let
$h(z)=re^{i(\varphi +\alpha (r)\pi )}$

For the set $\{z \in \C \ | \ 2 \leq |z|\}$,  $h(z)={\rm Id}$,
and for the set $\{z \in \C \ | \ |z|\leq \frac{3}{2}\}$, $h(z)$ a rotation by 
$180 ^{\circ}$ in the positive direction.

The diffeomorphism $h$ defined above induces an automorphism on $\pi _1(D \setminus K,u)$,
that switches the position of two generators of $\pi _1(D \setminus K,u)$, as
can be seen in the figure: 

\begin{center}
\epsfysize=4cm
\epsfbox{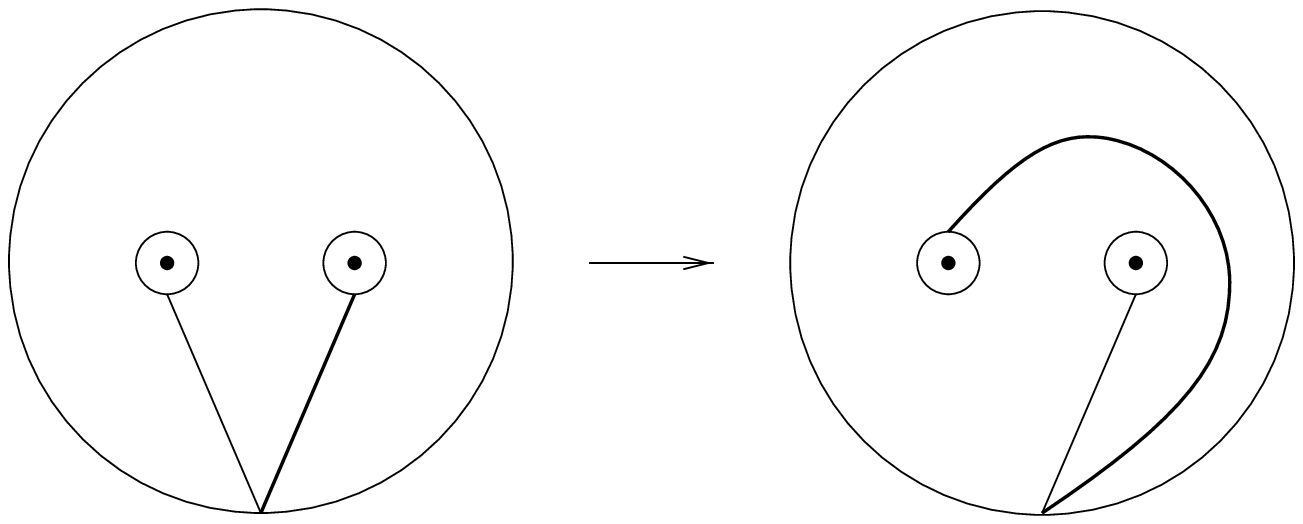}
\end{center}

Considering $(f \circ h \circ f^{-1})|_D$ (we will compose from left to
right) we get a diffeomorphism of $D$ which switches $a$ and $b$ and is the identity on $D
\setminus U$. Thus it defines an element of $B_n[D,K]$.

\bde
Let $H(\sigma)$ be the braid defined by $(f \circ h \circ f^{-1})|_D$. We call
$H(\sigma )$ the \underline{positive half-twist defined by $\sigma$}.
\ede

The half-twists generate $B_n$. In fact, one can choose $n-1$ half-twists that generates $B_n$ (see below):

\bde
Let $K=\{k_1,...,k_n\}$, and $\sigma _1,...,\sigma _{n-1}$ be a system of simple paths in
$D \setminus \partial D$ such that each $\sigma _i$ connects $k_i$ with $k_{i+1}$ and 

for all $i,j \in \{1,...,n-1\}$, $i<j$, $\left \{ 
\matrix{\sigma _i \cap \sigma _j=\emptyset & 2 \leq |i-j| \cr
\sigma _i \cap \sigma _{i+1}=\{k_{i+1}\} & i=1,...,n-2} \right.$. 

Let $H_i=H(\sigma _i)$. The ordered system of
(positive) half twists $(H_1,...,H_{n-1})$ are called a \underline{frame of $B_n[D,K]$}.
\ede

\bth
If $(H_1,...,H_{n-1})$ is a frame of $B_n[D,K]$, then
$B_n[D,K]$ is generated by $\{H_i\}_{i=1}^{n-1}$. Moreover, 
if $(H_1,...,H_{n-1})$ is a frame of $B_n[D,K]$, then the set
$\{H_i\}_{i=1}^{n-1}$ with the two relations $H_iH_j=H_jH_i$ if $2 \leq |i-j|$ and
$H_iH_{i+1}H_i=H_{i+1}H_iH_{i+1}$ for any $i=1,...,n-2$ are sufficient enough to present 
$B_n[D,K]$ and therefore this definition and Artin's definition for the braid group are equivalent.
\eth

\noindent
{\it Proof:} See \cite{BGTI}.

\bde
Let $\sigma _i$ be the straight line connecting $k_i$ to $k_{i+1}$, then we call the frame
$(H_1,...,H_{n-1})$ the \underline{standard frame} of $B_n[D,K]$.
\ede

Let $G=<\gamma _1,...,\gamma _n>$ be the ordered set of elements of the fundamental group that consists of the elements
of the g-base. Now, we look at the action of the diffeomorphisms in $B$ on $G$.
We can decide whether two diffeomorphisms $\beta _1$ and $\beta _2$ represent the same braid or not 
by looking at their action on $G$, since the automorphism on $\pi_1(D \setminus K,u)$ is
determined  by the action on its generators.

Because the set of half-twists in every frame generate $B_n[D,K]$, it is sufficient to check the 
action of each half-twist on $G$.

\bpr
Let $(H_1,...,H_{n-1})$ be the standard frame of $B_n[D,K]$ and $G=<\gamma _1,...,\gamma _n>$ be the 
ordered set of elements of the g-base for $\pi _1(D \setminus K,u)$. The action of $H_i$ on
$G$ is given by $G'=<\gamma _1,...,\gamma _{i-1},\gamma _{i+1},\gamma _{i+1}\gamma _i\gamma_{i+1}^{-1},\gamma
_{i+2},...,\gamma _n>$
\epr

\noindent
{\it Proof:} 
Follows immediately from the definition of the half-twist.
\hfill $\qed$

The last proposition means that when we work with the g-base elements we actually working
with a subset of the fundamental group, which is the conjugacy classes of the elements of
the standard g-base. We will exploit this characteristic in the future when we will present
the solution to the word problem.

\subsection{The braid word problem}
\bde
Let $b \in B_n$ be a braid. Then it is clear that $b=\sigma _{i_1}^{e_1} \cdot ...\cdot \sigma
_{i_l}^{e_l}$ for some sequence of generators, where $i_1,...,i_l \in \{1,...,n-1\}$ and 
$e_1,...,e_l \in \{1,-1\}$. We will call such a presentation of $b$ a \underline{braid
word}, and $\sigma _{i_k}^{e_k}$ will be called the \underline{$k^{th}$ letter of the word $b$}. 
$l$ is the \underline{length} of the braid word.
\ede

We will distinguish between two relations on the braid words.

\bde
Let $w_1$ and $w_2$ be two braid words. We will say that $w_1=w_2$ if they represent the
same element of the braid group.
\ede

\bde
Let $w_1$ and $w_2$ be two braid words. We will say that $w_1 \equiv w_2$ if $w_1$ and
$w_2$ are identical letter by letter.
\ede

Now, we can introduce the word problem: Given two braid words $w_1$ and $w_2$, 
decide whether $w_1=w_2$ or not.

\section{The computerized representation of the g-base}

In this section, we will describe the way we encode the g-base in \cite{WORD1}. It involves some
conventions.

Let $D$ be the closed unit disk, the point $u$ is the point $(0,-1)$ and the points
in $K$ are on the $x$-axis.

In order to encode the path in $D$, which is an element of the g-base, we will
distinguish some positions in $D$.

\bnot
We will denote by $(i,1)$ a point {\bf close} to $k_i$ but above it, $(i,-1)$ a point {\bf close} to 
$k_i$ but below it, and $(i,0)$ the point $k_i$ itself.
We will also denote the point $u$ by $(-1,0)$ (which is not its position in $D$, rather 
only a notation).
\enot

To represent a path in $D$, we will use a linked list which its links are based on the notations above, which
represents the position of the path in relation to the points $u$ and $k_i$, $i=1,...,n$.

Each link of the list holds the two numbers as described above. We will call them
(point,position).

\bex
The list $(1,0) \to (2,1) \to (3,1) \to (4,-1) \to (5,0)$ represents the following path:
\eex

\begin{center}
\epsfysize=1cm
\epsfbox{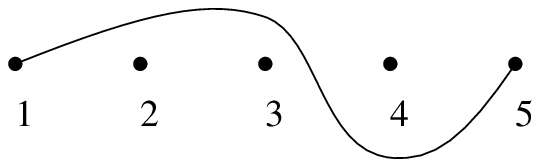}
\end{center}

As a rule, we will never connect the point $u$ to any point $(i,-1)$. This will
be done in order to obtain a unique way of representation, and to make the automatic 
computation of the twists easier.

We will be able to tell whether a path $(-1,0) \to (i,1)$ is passing to the left or to the 
right of the point $i$ simply by checking its continuation. 
If the path is turning to the left $(-1,0) \to (i,1) \to (i-1,e)$, then it is passing to the right
of the point $i$, and if the path is turning to the right $(-1,0) \to (i,1) \to (i+1,e)$, then it
is passing to the left of the point $i$ (where $e \in \{-1,1,0\}$).

\bex
The list $(-1,0) \to (3,1) \to (2,0)$ represents the following path:

\begin{center}
\epsfysize=3cm
\epsfbox{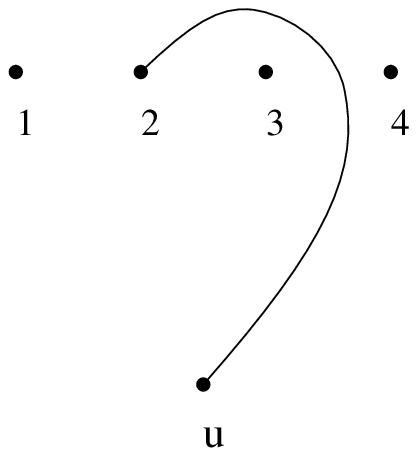}
\end{center}

The list $(-1,0) \to (2,1) \to (3,1) \to (3,-1) \to (2,0)$ represents the following path:

\begin{center}
\epsfysize=3cm
\epsfbox{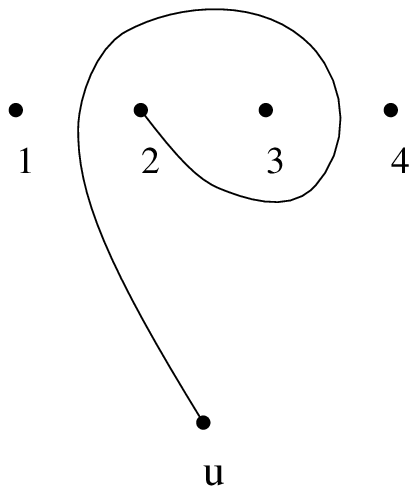}
\end{center}

\eex

In order to unify our treatment on all the paths of the g-base, we concatenate all 
of them into one list, which means that after we arrive at the end of one path (i.e. a
link $(i,0)$), the following link will be $(-1,0)$ marking the beginning of the next path.
For convenience, and not for mathematical reasons, we add the link $(-1,0)$ at the end of the
list.

\bex
The list $(-1,0) \to (1,1) \to (2,0) \to (-1,0) \to (1,0) \to (-1,0) \to (4,0) \to (-1,0) 
\to (4,1) \to (3,0) \to (-1,0)$ represents the g-base in the following figure 
(the small circles around the points are omitted):

\begin{center}
\epsfysize=3cm
\epsfbox{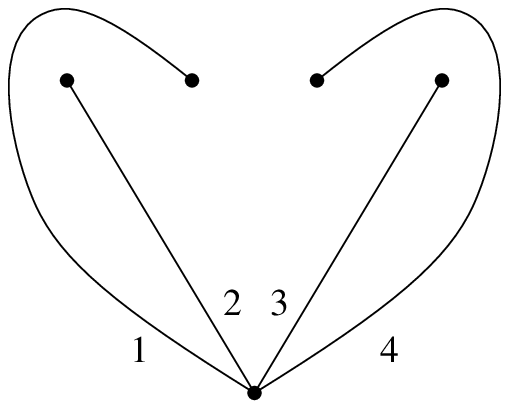}
\end{center}

\eex

The solution for the braid word problem given in \cite{WORD1} is based on an algorithm
that performs the action of the half-twists on the elements of the g-base, maintaining a
reduced form which is unique. If at the end of the process we find that the two g-bases
after the action of the two braid words are the same, then the two elements of the braid
group represent the same braid, and vice versa. This is without any questions, a
geometrical solution, and in the following sections we will present the syntactical
solution, using the same methods.

We will lay out an algorithmic way which will make it possible for us to apply the
group action of the braid group on the g-bases, thus obtaining a new presentation of the
braid group as an ordered set of elements of a g-base with the new operation.

In order to achieve this, we need a method that will enable us to write the elements of the g-base after the
action using the elements of the standard g-base, and another method that will make it
possible to return from the syntactic writing using the elements of the standard g-base
into the representation of its paths. This is going to be the main subject of the next two
sections.

\section{Transforming the geometric presentation to a syntactic one}

Let the standard g-base be $<\gamma _1,...,\gamma _n>$. In this section we are going to introduce 
the algorithms that gets as an input a path $\gamma \in \pi _1(D \setminus K,u)$ 
which is an element of the standard g-base after the action of some braid,
and returns a word in $\{\gamma _1,...,\gamma _n,\gamma _1^{-1},...,\gamma
_n^{-1}\}$ that represents the same element $\gamma$. 

We need more to say that every link $L$ consists of two numbers (Point,Position) and that
by writing $L.Point$ we mean the Point element in the link $L$, and by writing
$L.Position$ we mean the Position element of the link $L$.

\subsection{The $PathToSyntactic(\gamma)$ algorithm}
\balg
\parindent0pt
\parskip0pt

PathToSyntactic($\gamma$)

{\bf input}: $\gamma$ - a list that represents an element in $\pi _1(D \setminus K,u)$
which is an element of the standard g-base after the action of a braid.

{\bf output}: a word $w$ that consists of the generators of the fundamental group from the
standard g-base, which is the same element as $\gamma$.

{\bf PathToSyntactic($\gamma$)}

$w \leftarrow Null$

{\bf for} all the links $L$ in $\gamma$ {\bf do}

$\qquad$ $LastLink \leftarrow$ the link before $L$

$\qquad$ $NextLink \leftarrow$ the link after $L$

$\qquad$ $SecondLink \leftarrow$ the link after $NextLink$

$\qquad$ {\bf if} $L.Position=0$ {\bf then}

$\qquad$ $\qquad$ concatenate $\gamma _{L.Point}$ to $w$

$\qquad$ $\qquad$ {\bf continue}

$\qquad$ {\bf if} $L.Position=-1$ {\bf then}

$\qquad$ $\qquad$ {\bf continue}

$\qquad$ {\bf if} $L.Position=1$ {\bf then}

$\qquad$ $\qquad$ {\bf if} $NextLink.Point=L.Point-1$ {\bf then}

$\qquad$ $\qquad$ $\qquad$ concatenate $\gamma _{L.Point}$ to $w$

$\qquad$ $\qquad$ $\qquad$ {\bf continue}

$\qquad$ $\qquad$ {\bf if} $NextLink.Point=L.Point+1$ {\bf then}

$\qquad$ $\qquad$ $\qquad$ concatenate $\gamma _{L.Point}^{-1}$ to $w$

$\qquad$ $\qquad$ $\qquad$ {\bf continue}

$\qquad$ $\qquad$ {\bf if} $NextLink.Point=L.Point$ {\bf then}

$\qquad$ $\qquad$ $\qquad$ {\bf if} $SecondLink.Point=L.Point-1$ {\bf then}

$\qquad$ $\qquad$ $\qquad$ $\qquad$ concatenate $\gamma _{L.Point}^{-1}$ to $w$

$\qquad$ $\qquad$ $\qquad$ $\qquad$ {\bf continue}

$\qquad$ $\qquad$ $\qquad$ {\bf if} $SecondLink.Point=L.Point+1$ {\bf then}

$\qquad$ $\qquad$ $\qquad$ $\qquad$ concatenate $\gamma _{L.Point}$ to $w$

$\qquad$ $\qquad$ $\qquad$ $\qquad$ {\bf continue}

$l \leftarrow$ the length of $w$

{\bf for} all $i$ {\bf from} $l-1$ {\bf to} $1$ {\bf step} $-1$ {\bf do}

$\qquad$ $t \leftarrow$ the $i^{th}$ letter in $w$

$\qquad$ concatenate to $w$ the letter $t^{-1}$

{\bf return} $w$

\ealg
\parindent10pt
\parskip10pt

\subsection{Proof of correctness}

\bpr
Let $\gamma$ be the list representing the path. Then, any link in $\gamma$ that has position 
$-1$ does not contribute any letter to the word $w$ representing the same elements 
using the standard g-base elements as generators.
\epr

\noindent
{\it Proof:}
Let us look at a part of the path that goes beneath the points $(i,0)$. This part consists only on 
links of the type $(i,-1)$. The link that is before this sublist must be of the type $(j,1)$, because 
it can not be the link denoted $(-1,0)$, which is the first link in the path, 
and must be followed by convention with a link above a point, 
nor any of the links $(j,0)$ since those links represents the end 
of the path.

We look at the path that this part represents. This is a line going from the one point (say 
$i$) to another (say $j$), which goes below the points of $K$, and therefore below 
the x-axis. This means that the part of the path represented by those links is homotopically 
equivalent to two straight lines forming a 'V' shape starting from $(i,1) \to (-1,0) \to (j,e)$ where
$e \in \{0,1\}$. These two lines are represented by two elements of the standard g-base 
that are added by looking at links just before this section and immediately after it, hence, this part of the path does not 
contribute any letters to the word $w$.
\hfill $\qed$

\bpr
Let $\gamma$ be the list representing the path. The last link $L=(i,0)$ of $\gamma$ contributes 
the letter $\gamma _i$ to $w$.
\epr

\noindent
{\it Proof:}
The last link in the list representing a path of the g-base is always of the type $(i,0)$ 
It represents the loop of the path around the point $i$. This contributes the element of the standard 
g-base that corresponds to that loop, which is $\gamma _i$.
\hfill $\qed$

\bpr
Let $L=(i,1)$ be a link in the list $\gamma$. The element of the standard g-base 
that corresponds to $L$ and therefore we need to add to the word $w$ is determined by 
the following rules:

\ben
\item
If the path goes to the left after the link. Then, we have to add the letter $\gamma _i$.

\item
If the path goes to the right after the link. Then, we have to add the letter $\gamma _i^{-1}$.

\item
If the path goes to the same point then, the letter we have to add is determined by the next link.
If the following link goes to the left we have to add the letter $\gamma _i^{-1}$, and if the 
following link goes to the right we have to add the letter $\gamma _i$.
\een
\epr

\noindent
{\it Proof:}
The easiest way to prove this is to look at the homotopy type geometrically.
In case 1, the path simply goes from right to left above the point $i$ which means that homotopically 
all we need to add is $\gamma _i$, (see the following figure (a)).

Case 2, is the same as case 1 but on the other direction, which means that homotopically 
we need to add $\gamma _i^{-1}$, (see the following figure (b)).

Case 3, is somewhat more complicated. The only reason why we might have two consecutive links 
over the same point is when we are switching direction around this point. Therefore if the path 
goes to the left, it came from the left, and rounded the point $i$ clockwise (because the 
initial point was $(i,1)$). This means that we must add the letter $\gamma _i^{-1}$. 
(see the following figure (c)).

On the other hand, if the path continues to the right, it came from the right, 
and rounded the point $i$ counterclockwise (again because the initial point was above $i$). 
This means that we have to add the letter $\gamma _i$. (see the following figure (d)).
\hfill $\qed$

\begin{center}
\epsfysize=5.5cm
\epsfbox{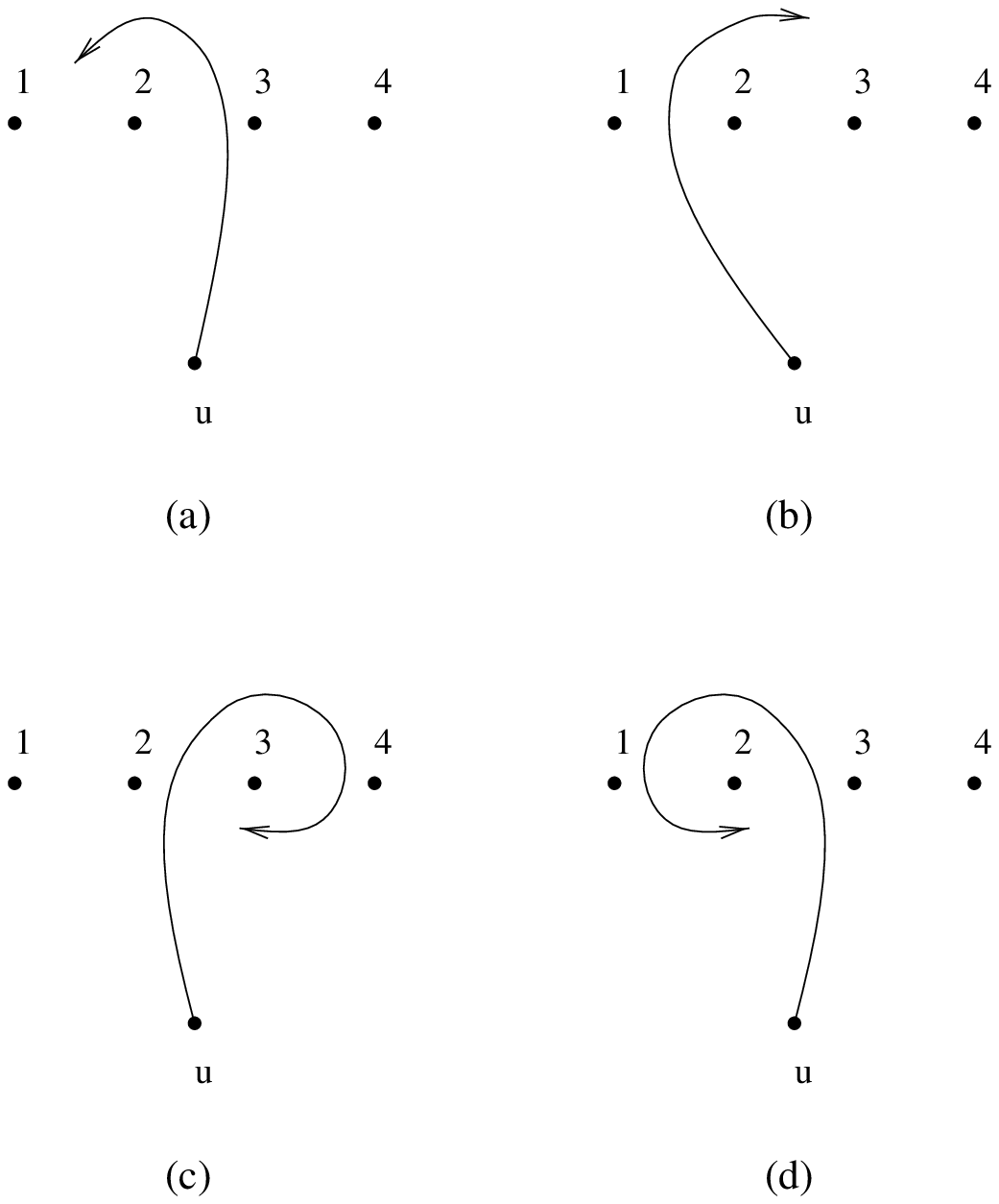}
\end{center}

\bpr
Let $w$ be the word resulted by going over the list and $l$ be the length of $w$. 
Let $w'$ be the word created from $w$ by taking the first $l-1$ letters of $w$.  
Then, in order to complete the word $w$ one has to make the word to a conjugacy by concatenating 
$w'^{-1}$ to $w$.
\epr

\noindent
{\it Proof:}
This is true because the path is encoded in such a way that the way back after the circle around 
the end point is omitted. But, the way back is exactly the opposite of the way to the end point, 
and therefore can be represented exactly by the negation of the same word.
\hfill $\qed$

\bth
The algorithm $PathToSyntactic(\gamma)$ returns a word in the elements of the standard g-base 
that represents the same element of the fundamental group represented by $\gamma$.
\eth

\hfill $\qed$

This concludes the proof of the correctness of the algorithm. In the next section we will 
compute its complexity.

\subsection{Complexity}
In this subsection we will compute and prove the complexity and the amount of memory
needed for the $PathToSyntactic(\gamma)$ algorithm.

\bpr
Let $n$ be the length of the list representing the path $\gamma$. Then, 
the complexity of the algorithm $PathToSyntactic(\gamma)$ is bounded by $O(n)$ operations.
\epr

\noindent
{\it Proof:}
The algorithm goes over every link in the list exactly once. On each of the links it may at most 
add one element of the standard g-base to the word $w$. In the worst case this yields $n$ operations 
exactly. Now, at the end of this process the algorithm doubles the result word $w$ to make it conjugated. 
The total of operations therefore, is bounded by $2n$ which is $O(n)$.
\hfill $\qed$

\bpr
Let $n$ be the length of the list representing the path $\gamma$. Then, 
The amount of memory needed for the algorithm $PathToSyntactic(\gamma)$ is bounded by $O(n)$.
\epr

\noindent
{\it Proof:}
Apart from a fixed number of counters and help variables, the algorithm $PathToSyntactic(\gamma)$ 
uses only the memory needed to store the word $w$, which is linear in $n$. Therefore the 
amount of memory it needs is bounded by $O(n)$.
\hfill $\qed$

\section{Transforming the syntactic presentation to the geometric one}

Let the standard g-base be $<\gamma _1,...,\gamma _n>$. In this section we are going to introduce 
the algorithm that gets as an input a word in the generators of the fundamental group $\{\gamma _1
,...,\gamma _n\}$, and will return a list that represents the path in $\pi _1(D \setminus
K,u)$. All paths will be encoded in the way that has been described in section 2.

\subsection{The SyntacticToPath($w$) algorithm}
\balg
\parindent0pt
\parskip0pt

SyntacticToPath($w$)

{\bf input}: a word $w$ that consists of the generators of the fundamental group from the
standard g-base, which is a conjugacy of an element of $\gamma$, written such that $w
\equiv Q^{-1}\gamma Q$.

{\bf output}: a list that represents an element in $\pi _1(D \setminus K,u)$
which is an element of the standard g-base after the action of a braid, which is the same
as $w$.

{\bf SyntacticToPath($w$)}

{\bf if} $w=Null$ {\bf then}

$\qquad$ {\bf return} the standard g-base

$L \leftarrow$ the empty list

set first link in $L$ to $(-1,0)$

{\bf for} every letter $l$ in the first $\frac{|w|-1}{2}$ letters of $w$ {\bf do}

$\qquad$ $CurrentPoint \leftarrow$ the generator index in $l$

$\qquad$ $NextPoint \leftarrow$ the generator index in the letter after $l$

$\qquad$ $LastLetter \leftarrow$ the letter before $l$

$\qquad$ {\bf if} $l$ is the first letter in $w$ {\bf then}

$\qquad$ $\qquad$ add to $L$ the link $(CurrentPoint,1)$

$\qquad$ $\qquad$ $LastPoint \leftarrow CurrentPoint$

$\qquad$ $\qquad$ {\bf continue}

$\qquad$ {\bf if} $LastLetter=\gamma _i^{-1}$ for some $i$ {\bf then}

$\qquad$ $\qquad$ {\bf if} $LastPoint<CurrentPoint$ {\bf then} 

$\qquad$ $\qquad$ $\qquad$ {\bf for} all $c$ {\bf from} $LastPoint+1$ {\bf to} $CurrentPoint-1$ {\bf
do}

$\qquad$ $\qquad$ $\qquad$ $\qquad$ add to $L$ the link $(c,-1)$

$\qquad$ $\qquad$ $\qquad$ {\bf if} $l=\gamma _i$ for some $i$ {\bf then}

$\qquad$ $\qquad$ $\qquad$ $\qquad$ add to $L$ the link $(CurrentPoint,-1)$

$\qquad$ $\qquad$ $\qquad$ add to $L$ the link $(CurrentPoint,1)$

$\qquad$ $\qquad$ $\qquad$ {\bf continue}

$\qquad$ $\qquad$ {\bf else}

$\qquad$ $\qquad$ $\qquad$ {\bf for} all $c$ {\bf from} $LastPoint$ {\bf to} $CurrentPoint+1$
{\bf step} $-1$ {\bf do}

$\qquad$ $\qquad$ $\qquad$ $\qquad$ add to $L$ the link $(c,-1)$

$\qquad$ $\qquad$ $\qquad$ {\bf if} $l=\gamma _i^{-1}$ for some $i$ {\bf then}

$\qquad$ $\qquad$ $\qquad$ $\qquad$ add to $L$ the link $(CurrentPoint,-1)$

$\qquad$ $\qquad$ $\qquad$ add to $L$ the link $(CurrentPoint,1)$

$\qquad$ $\qquad$ $\qquad$ {\bf continue}

$\qquad$ {\bf else}

$\qquad$ $\qquad$ {\bf if} $LastPoint>CurrentPoint$ {\bf then}

$\qquad$ $\qquad$ $\qquad$ {\bf for} all $c$ {\bf from} $LastPoint-1$ {\bf to} $CurrentPoint+1$
{\bf step} $-1$ {\bf do}

$\qquad$ $\qquad$ $\qquad$ $\qquad$ add to $L$ the link $(c,-1)$

$\qquad$ $\qquad$ $\qquad$ {\bf if} $l=\gamma _i^{-1}$ for some $i$ {\bf then}

$\qquad$ $\qquad$ $\qquad$ $\qquad$ add to $L$ the link $(CurrentPoint,-1)$

$\qquad$ $\qquad$ $\qquad$ add to $L$ the link $(CurrentPoint,1)$

$\qquad$ $\qquad$ $\qquad$ {\bf continue}

$\qquad$ $\qquad$ {\bf else}

$\qquad$ $\qquad$ $\qquad$ {\bf for} all $c$ {\bf from} $LastPoint$ {\bf to} $CurrentPoint-1$
{\bf do}

$\qquad$ $\qquad$ $\qquad$ $\qquad$ add to $L$ the link $(c,-1)$

$\qquad$ $\qquad$ $\qquad$ {\bf if} $l=\gamma _i$ for some $i$ {\bf then}

$\qquad$ $\qquad$ $\qquad$ $\qquad$ add to $L$ the link $(CurrentPoint,-1)$

$\qquad$ $\qquad$ $\qquad$ add to $L$ the link $(CurrentPoint,1)$

$\qquad$ $\qquad$ $\qquad$ {\bf continue}

set the last link in $L$ Position to 0

add to $L$ the link $(-1,0)$

Reduce(L)

{\bf return} $L$

\ealg
\parindent10pt
\parskip10pt

\subsection{Proof of correctness}
In this subsection we will prove that the algorithm $SyntacticToPath(w)$ returns a list
representing an element of a g-base which is homotopically equivalent to the element
represented by the word $w$ that consists of the generators of the standard g-base.

\bpr
Let $w$ be the input word consists of the generators of the fundamental group. Then, in order to compute 
the path represented by the word $w$ it is sufficient enough to look at the first half of
$w$.
\epr

\noindent
{\it Proof:}
The word $w$ represents a conjugated element in the fundamental group of one of the elements
of the standard g-base, and written as $w \equiv Q^{-1}\gamma _iQ$ for some element of the fundamental group $Q$.
Therefore, the second half of the word $w$ represents exactly the back path of the loop
represented by $w$. Since by convention we represent the path only by a sequence of links
that start at the base point $u$ and ends at one of the points of $K$, the back path is
omitted and therefore, not necessary for the computation. 
\hfill $\qed$

\bpr
Let $w$ be the input word consists of the generators of the fundamental group. If the
first letter in $w$ is $\gamma _i$ then, the second link in the list must be $(i,1)$ or
$(i,0)$.
\epr

\noindent
{\it Proof:}
The element $\gamma _i$ represents a loop that rounds the point $i$. $w=\gamma _i$ then,
the list the algorithm has to return is $(-1,0) \to (i,0) \to (-1,0)$. So, we need to
check only the case when there are more generators in $w$. In that case, since all the
generators in $w$ connect always at $u$, This means that the loop $\gamma _i$ represents
the fact that the path goes over the point $i$ and into one of the directions left or right. 
\hfill $\qed$

\bpr
Let $\gamma _i$ and $\gamma _j$ be two consecutive letters in $w$, and suppose $i<j$. 
Then, one of the following cases must happen:
\ben

\item
If $\gamma _i$ is in a negative power, and $\gamma _j$ is in a positive power. 
Then, we have to add the links $(i+1,-1) \to ... \to (j,-1) \to (j,1)$.

\item
If both $\gamma _i$ and $\gamma _j$ are in a negative power. 
Then, we have to add the links $(i+1,-1) \to ... \to (j-1,-1) \to (j,1)$.

\item
If $\gamma _i$ is in a positive power, and $\gamma _j$ is in a negative power. 
Then, we have to add the links $(i,-1) \to ... \to (j-1,-1) \to (j,1)$.
\item

If both $\gamma _i$ and $\gamma _j$ are in a positive power. 
Then, we have to add the links $(i,-1) \to ... \to (j,-1) \to (j,1)$.
\een
\epr

\noindent
{\it Proof:}
The proof follows from the homotopy relation between each of the cases and the list of the
links we add.
\ben
\item
In this case we have a negative loop over the point $i$, which lies to the left of the
point $j$, which is rounded by the second loop. This means that we orbit the point $j$
in a counterclockwise direction, after coming from the left. This is homotopic to the
list of links $(i+1,-1) \to ... \to (j,-1) \to (j,1)$, as can be seen by the following
figure:

\begin{center}
\epsfysize=2cm
\epsfbox{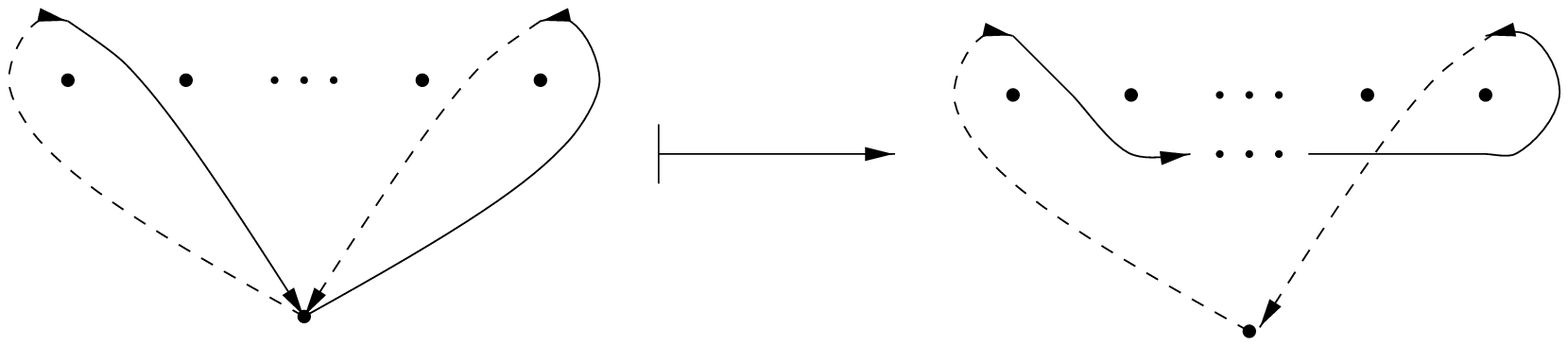}
\end{center}

\item
In this case we have a negative loop over the point $i$, which lies to the left of the
point $j$, which is rounded by the second negative loop. This means that we simply have a
path going from left to right below the points in between $i$ and $j$. This is homotopic
to the list of links $(i+1,-1) \to ... \to (j-1,-1) \to (j,1)$, as can be seen by the
following figure:

\begin{center}
\epsfysize=2cm
\epsfbox{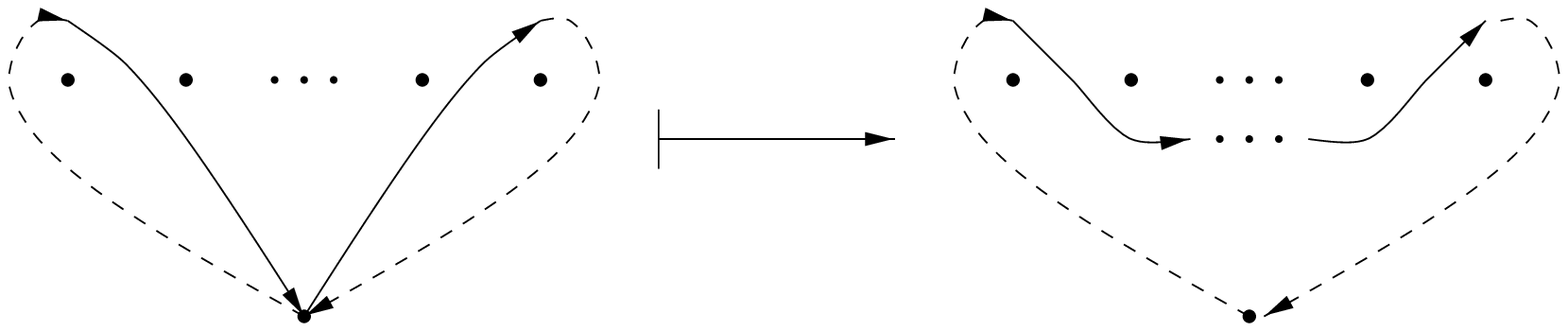}
\end{center}

\item
In this case we have a positive loop over the point $i$, which lies to the left of the
point $j$, which is rounded by the second negative loop. This means that we orbit the
point $i$ counterclockwise and then, we go below the points between $i$ and $j$. this is
homotopic to the list of links $(i,-1) \to ... \to (j-1,-1) \to (j,1)$, as can be seen
in the following figure:

\begin{center}
\epsfysize=2cm
\epsfbox{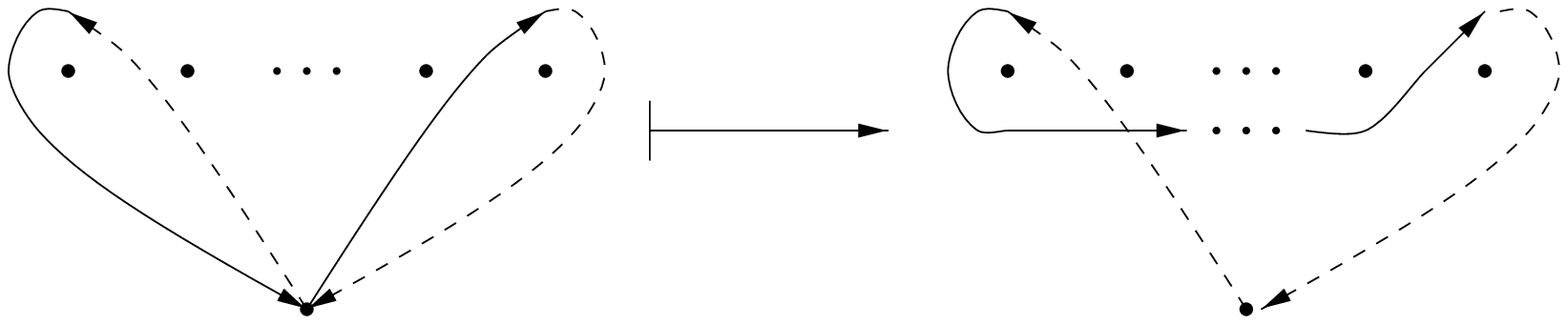}
\end{center}

\item
In this case we have a positive loop over the point $i$, which lies to the left of the
point $j$, which is rounded by the second positive loop. This means that we that we orbit
both the point $i$ and the point $j$ in a counterclockwise direction. This is homotopic to
the list of links $(i,-1) \to ... \to (j,-1) \to (j,1)$, as can be seen by the following
figure:

\begin{center}
\epsfysize=2cm
\epsfbox{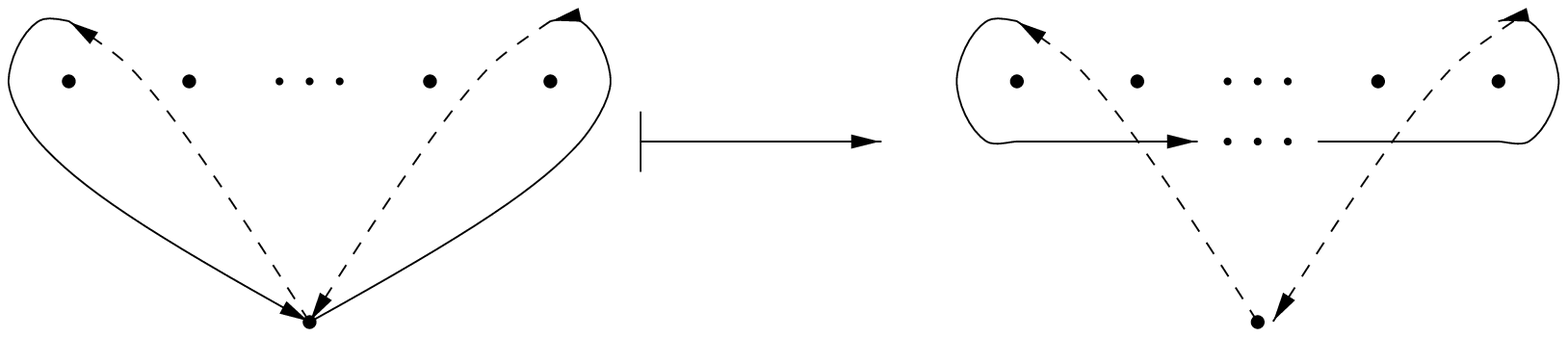}
\end{center}

\een
\hfill $\qed$

Now, we will consider the case when $j<i$, which means that the path is going from left to
right, so we have the following proposition:

\bpr
Let $\gamma _i$ and $\gamma _j$ be two consecutive letters in $w$, and suppose $j<i$. 
Then, one of the following cases must happen:
\ben

\item
If $\gamma _i$ is in a negative power, and $\gamma _j$ is in a positive power. 
Then, we have to add the links $(i,-1) \to ... \to (j+1,-1) \to (j,1)$.

\item
If both $\gamma _i$ and $\gamma _j$ are in a negative power. 
Then, we have to add the links $(i,-1) \to ... \to (j,-1) \to (j,1)$.

\item
If $\gamma _i$ is in a positive power, and $\gamma _j$ is in a negative power. 
Then, we have to add the links $(i-1,-1) \to ... \to (j,-1) \to (j,1)$.
\item

If both $\gamma _i$ and $\gamma _j$ are in a positive power. 
Then, we have to add the links $(i-1,-1) \to ... \to (j+1,-1) \to (j,1)$.
\een
\epr

\noindent
{\it Proof:}
The proof is similar to the proof of the proposition $4.4$ above.
\hfill $\qed$

Now, what is left is to notice that the middle letter in $w$ marks the end point of the
path. This is due to the convention of the representation of the path that marks the end
of the path in a link connected to the point $i$, that is $(i,0)$. So, what the algorithm
does is when finishing passing over the letters, it simply connects the last link to the
point it is associated with. In order to comply with the representation of the path's
conventions, the algorithm adds in the end one more link $(-1,0)$, and reduces its
presentation by using the $Reduce(L)$ algorithm found in \cite{WORD1}.

This completes the proof of the correctness of the algorithm, so we get:

\bth
The algorithm $SyntacticToPath(w)$ returns a list representing the the same element of the fundamental 
group as the word $w$ in the elements of the standard g-base represents, where $w \equiv
Q^{-1}\gamma Q$.
\eth
\hfill $\qed$

\subsection{Complexity}
In this subsection we will compute and prove the complexity and the amount of memory
needed for the $SyntacticToPath(w)$ algorithm.

\bpr
Let $l$ be the length of the word $w$ in the standard set of generators of the fundamental
group, and $n$ be the number of points in $K$. Then, the complexity of the algorithm $SyntacticToPath(w)$ 
is bounded by $O(nl)$ operations.
\epr

\noindent
{\it Proof:}
The algorithm goes over every letter in the first half of the word $w$ exactly once. By
each letter the algorithm might add at most $n+1$ links to the list, since it might add all the
links between $i$ and $j$ where $\gamma _i$ and $\gamma _j$ are consecutive letters in
$w$, and since $|i-j| < n+1$. At the end of the algorithm it uses the $Reduce(g)$ algorithm 
that was presented in \cite{WORD1}. Its complexity is bounded by the length of the list as proven in 
\cite{WORD1}. So, at most the algorithm complexity is bounded by $nl$ operations which is $O(nl)$.
\hfill $\qed$

\bpr
Let $l$ be the length of the word $w$ in the standard set of generators of the fundamental
group, and $n$ be the number of points in $K$. Then, 
The amount of memory needed for the algorithm $SyntacticToPath(w)$ is bounded by $O(nl)$.
\epr

\noindent
{\it Proof:}
Besides a finite predetermined number of variables, the algorithm creates a new linked
list, that represents the path. The size of this list was calculated in the proof of the
complexity proposition. Therefore, the list might in the worst case have $(n+1)l+2$ links,
which is bounded by $O(nl)$.
\hfill $\qed$

\section{Group structure over the presentation of the g-base}
When we consider the lists that represent elements of the g-base, as described in
\cite{WORD1}, we get a unique way for presenting elements of the braid group. So, 
we thought of a way to multiply lists representing different g-bases. By multiplying, or
by adding a group structure to the presentation of the g-bases, we mean that when we have
two lists $L_1,L_2$ representing braids $\beta_1,\beta_2$, (that is, two lists that each one 
of them represents a g-base, which is the result of the action of the braid on the
standard g-base), we want to be able to perform a process that will result with a list
$L_3=L_1L_2$ that represents a g-base, which is the result of the action of the braid
$\beta_1 \beta_2$ on the standard g-base. 

Since multiplying braids is actually the
composition of diffeomorphisms, the result on the g-base is a composition of the elements
of the g-base. When we have two lists $L_1$ and $L_2$ that represents the two braids
$\beta _1$ and $\beta _2$, we have to look at the action of the first braid
on the standard g-base $<\gamma _1,...,\gamma _n>$, resulting with a different g-base 
say $<\gamma _1',...,\gamma _n'>$ and than to act with $\beta _2$ on the elements of the
new g-base, resulting with a third g-base say $<\gamma _1'',...,\gamma _n''>$. The latter is
the result of the action of the braid $\beta _1 \beta _2$ on the standard g-base.

This is rather difficult to perform from a geometrical point of view, but fortunately 
easy to perform on the syntactical presentation of the elements of the g-base. The action
is as follows: 

\balg
\parindent0pt
\parskip0pt

Multiply($L_1,L_2$)

{\bf input}: $L_1$ and $L_2$ - two lists representing two g-bases associated with $\beta
_1$ and $\beta _2$.

{\bf output}: A list representing the g-base resulted after the action of the braid $\beta
_1 \beta _2$ on the standard g-base.

{\bf Multiply($L_1,L_2$)}

Let $n$ be the number of points in $K$

Break $L_1$ into $n$ sublists $L_1[i]$ each one representing one element of the g-base.

Break $L_2$ into $n$ sublists $L_2[i]$ each one representing one element of the g-base.

{\bf for} all $i$ {\bf from} 1 {\bf to} $n$ {\bf do}

$\qquad$ $w_1[i] \leftarrow PathToSyntactic(L_1[i])$

$\qquad$ $w_2[i] \leftarrow PathToSyntactic(L_2[i])$

{\bf for} all $i$ {\bf from} 1 {\bf to} $n$ {\bf do}

$\qquad$ replace in $w_1[i] $every letter $\gamma _j$ with $w_2[j]$,

$\qquad$ and every letter $\gamma _j^{-1}$ with $w_2[j]^{-1}$

$\qquad$ Store the result in $w_3[i]$

$\qquad$ Eliminate every two consecutive letters in $w_3[i]$ of the type $\gamma
_i^e\gamma _i^{-e}$ where $e \in \{1,-1\}$.

$\qquad$ $L_3[i] \leftarrow SyntacticToPath(w_3[i])$

$L_3 \leftarrow$ the concatenation of all the paths in $L_3[i]$ into one list representing

the resulted g-base

{\bf return} $L_3$

\ealg
\parindent10pt
\parskip10pt

\bpr
Let $l_1$ and $l_2$ be the lengths of the lists $L_1$ and $L_2$, suppose $n$ is the number 
of strings in the braid group. Then, the complexity of the algorithm is bounded by $O(nl_1l_2)$
\epr

\noindent
{\it Proof:}
Breaking the lists into $n$ small lists takes $O(l_1+l_2)$ actions, since the sum of all
the lengths of all the paths $L_1[i]$ is equal to $l_1+n$ (We add the link
$(-1,0)$ to each path), and the sum of all the lengths of all the paths $L_2[i]$ is equal
to $l_2+n$, and  
the action of $PathToSyntactic(L_1[i])$ and $PathToSyntactic(L_2[i])$ is linear in the
length of all the lists.

The replacement of every letter in $w_1[i]$ with the word $w_2[j]$ is bounded by the
multiplication of the length of $w_1[i]$ by $l_2$. Again, since the sum of all the lengths
of $w_1[i]$ is $l_1+n$, and since we do not replace the elements $(-1,0)$, we result with a boundary of $O(l_1l_2)$.

The elimination process is linear in the length of $w_3[i]$ since each letter that was
inserted to $w_3[i]$ can be extracted only once. If one uses a doubly-connected list or a vector to
store the letters, then, moving in both direction on the list is possible resulting a
linear time elimination. Note that the size of all the elements $w_3[i]$ is bounded by the
previous argument by $O(l_1l_2)$.

Concatenation of all the paths takes $2n$ operations, since this is only the connection of
linked lists, or $O(l_1l_2)$ if they are stored in a vector.

The only thing we have not checked yet is how long it takes for the algorithm to transform every syntactic
word $w_3[i]$ into a list. But, as proved above this is bounded by $O(nl_1l_2)$.

So, at the end we have reached the boundary of $O(nl_1l_2)$.
\hfill $\qed$

\bpr
Let $l_1$ and $l_2$ be the lengths of the lists $L_1$ and $L_2$. Then, the amount of memory
needed for the algorithm is bounded by $O(l_1l_2)$.
\epr

\noindent
{\it Proof:}
The largest structure needed for the algorithm is the resulted linked list $L_3$ which
has $O(l_1l_2)$ links. All the other data structures are lists that the sum of their sizes
do not exceed $O(l_1 + l_2)$, and variables that their number is preknown and constant. Therefore, the proposition is proved.
\hfill $\qed$

\section{An improvement of the braid word solution}
In this section, we will show how one can improve the solution for the braid word problem
that was presented in \cite{WORD1}. The main idea of the new solution is to base it on the
syntactic presentation of the elements of the g-base. Since every element in the g-base, 
resulted by an action of a braid, is conjugated to an element of the standard g-base, it is
possible to encode the elements during the computation in a more efficient way.

Unfortunately, the process of performing the action of the diffeomorphisms on the elements of a g-base still takes a
long time. But, since the algorithm in \cite{WORD1} is very practical for short words, and
since using the new method will improve the running time of the algorithm, we result in an even more
practical solution to the word problem.

\bde
Let $\Gamma =<\gamma _1,...,\gamma _n>$ be a g-base. We say that the g-base $\Delta =<\gamma _1',...,\gamma
_n'>$ is obtained from $\Gamma$ by the \underline{braid move $H_i$} (or $\Gamma$ is
obtained from $\Delta$ by the braid move $H_i^{-1}$) if 

\ben
\item
$\gamma' _j=\gamma _j$ for all $j \neq i,i+1$
\item
$\gamma' _i=\gamma _{i+1}$
\item
$\gamma' _{i+1}=\gamma _{i+1}\gamma _i\gamma _{i+1}^{-1}$
\een
\ede

If we consider the elements of the standard frame $H_1,...,H_{n-1}$ as the generator set
for the braid group, then the action that $H_i$ performs over a g-base is exactly the
$i^{th}$ braid move on it.

Due to the fact that the group operation between two braids $\beta _1$ and $\beta _2$ is a 
composition of diffeomorphisms, in order to compute the action on the standard g-base (or any
g-base for that matter), we need to go over the braid word from the end to
the beginning and perform on the g-base at each step the braid move induced by the generator at this
position. This is the same as going over the braid word from the beginning to the end,
replacing for each letter $H_i$ each appearance of $\gamma _i$ with $\gamma _{i+1}$ and
$\gamma _{i+1}$ with $\gamma _{i+1} \gamma _i \gamma _{i+1}^{-1}$, or its negative compliant for $H_i^{-1}$.

This yields the following algorithm for encoding the diffeomorphisms action on the
standard g-base, using the braid move defined above.

\subsection{The ProcessWord(w) algorithm}

In this section, we will introduce the changed $ProcessWord(w)$ algorithm. This algorithm as
explained above will process on the syntactic presentation of the g-base. 

\balg
\parindent0pt
\parskip0pt

ProcessWord($w$)

{\bf input}: $w$ - a braid word.

{\bf output}: $w'$ - a word which represents the g-base elements written in the standard
generators of the fundamental group.

We denote $e \in \{1,-1\}$

{\bf ProcessWord($w$)}

$\Gamma \leftarrow$ the standard g-base $<\gamma _1,...,\gamma _n>$

{\bf for} each letter $l=H_i^e$ in the braid word $w$ {\bf from} the end {\bf to} the
beginning {\bf do}

$\qquad$ {\bf if} $e=1$ {\bf then}

$\qquad$ $\qquad$ activate the braid move $H_i$ on $\Gamma$

$\qquad$ {\bf else}

$\qquad$ $\qquad$ activate the braid move $H_i^{-1}$ on $\Gamma$

$\qquad$ $Reduce(\Gamma)$

{\bf return} $\Gamma$

\ealg
\parindent10pt
\parskip10pt

The action of the braid generator on the g-base, is given as above by maintaining $n-2$
elements without a change, replacing one elements position with the other, and
conjugating the latter. Therefore, we will show how to encode the conjugacy in a way that
will make the computation fast on one hand, and will save a lot of memory on the other.

\bnot
We will denote $B_i^j$ the part of the word $\gamma$ that represents an element of the
g-base, that begins at the $i^{th}$ letter and its length is $j$ letters.
\enot

This will make it possible to reduce the size of the memory needed to represents elements
of the changed g-base.

\bex
Instead of writing the element of the g-base $\gamma =\gamma _3\gamma _2\gamma _1\gamma
_2^{-1}\gamma _3^{-1}$ we can write $\gamma = \gamma _3\gamma _2\gamma _1B_1^2$
\eex

This method helps to keep the number of moves and copies we have to make small, 
although it makes it a bit more difficult to perform the $Reduce(\Gamma)$ procedure.

\subsection{The Reduce($\Gamma$) algorithm}

The new version of the $Reduce(\Gamma)$ algorithm gets as an input an ordered set of syntactic representations
of a g-base for the fundamental group, that we get after the operation of a braid generator
on it. The algorithm performs a reduction of these words into the smallest form possible.

Since the fundamental group of a punctured disk is a free group, the 
only reduction rule available is when we have to consecutive letters $\gamma _i\gamma _i^{-1}$
that we can immediately delete both.

So, what the $Reduce(\Gamma)$ algorithm does is to go over each element in $\Gamma$ and
eliminate all the pairs $\gamma _i \gamma _i^{-1}$ and $\gamma _i^{-1} \gamma _i$.

\bde
Let $\pi _1(D \setminus K,u)$ be the fundamental group of a punctured disk, and let $w$ be
a word in its standard generators. We call $w$ reduced if $w$ does not contain any two
consecutive letters of the form $\gamma _i^e \gamma _i^{-e}$ for any $1 \leq i \leq n$ and
$e \in \{-1,1\}$.
\ede

\subsection{Why is that better than before}
Since the fundamental group of a punctured disk is a free group on $n$ generators,
two reduced words $w$ and $w'$ are the same if and only if $w \equiv w'$ ($w$ and $w'$ are equal 
letter by letter).

This is easy to check and therefore, we hold a solution for the braid word problem. If one
wants to compare two braid words $\beta _1$ and $\beta _2$, he has to 
compute $\Gamma _1=ProcessWord(\beta _1)$ and $\Gamma _2=ProcessWord(\beta _2)$ and
compare the elements of $\Gamma _1$ with the elements of $\Gamma _2$.

Denote $\Gamma _1=<\gamma _{1,1},...,\gamma _{1,n}>$ and $\Gamma _2=<\gamma
_{2,1},...,\gamma _{2,n}>$, where $\gamma _{j,i}$ is reduced, than $\beta _1=\beta _2$ if and only if for every $1 \leq i
\leq n$ we have $\gamma _{1,i} \equiv \gamma _{2,i}$.

Although not in complexity, this is faster than the algorithm $ProcessWord(w)$ that was 
presented in \cite{WORD1}, since no maintaining of the linked list in necessary. Since the
algorithm there was very practical for short braid words, this one is even more practical.

\section{Conclusions}

Although for very long braid words this algorithm running time is long, due 
to the fact that the complexity of the g-base grows with the length of the braid word,
for short braid words we obtained a quick algorithm in comparison with other methods.

The worst drawback of this algorithm is the representation of the conjugated elements in
$\pi _1(D \setminus K,u)$ to the elements of the standard g-base. Perhaps a new method for
encoding these elements will make it possible to reduce the running time of the algorithm
even more.

Another thing, that can help to reduce the running time of the algorithm is to add
some generators to the standard generators of the fundamental group $\pi _1(D \setminus
K,u)$. These elements will make it possible to shorten the length of the words that
represents the different elements of the g-base. Some ideas we thought about are to add
the elements $\gamma _{i,j}=\gamma _i \cdot ... \cdot \gamma _j$. This set of generators
will make it possible to write an element of the g-base using only its path turning points,
which will shorten the word immensely.

We believe that there is a connection between presentations of two conjugated braid words,
therefore we believe that this connection might yield us a practical algorithm for solving
the conjugacy problem in the braid group.

The operation of the unprocess of the algorithm here, which means to compute the braid
word from a given g-base is even easier than in the case of the other presentation.

Using this new method allows us to leave the geometrical solution in \cite{WORD1}, this
makes things even easier in order to compute the braid monodromy automatically for some
of the cases.

\begin{\bib}{10}
\bibitem{Artin} Artin, E., {\it Theory of braids}, Ann. Math. {\bf 48} (1947), 101-126.
\bibitem{NEW} Birman, J.S., Ko, K.H. and Lee, S.J., {\it A new approach to the word and conjugacy problems in the braid groups}, Adv. Math. {\bf 139} (1998), 322-353.
\bibitem{Deh1} Dehornoy, P., {\it From large cardinals to braids via distributive algebra}, J. Knot Theory \& Ramifications {\bf 4(1)} (1995), 33-79.
\bibitem{Deh2} Dehornoy, P., {\it A fast method for comparing braids}, Adv. Math. {\bf 125(2)} (1997), 200--235. 
\bibitem{POSBR} Elrifai, E.A. and Morton, H.R., {\it Algorithms for positive braids}, Quart. J. Math. Oxford Ser. (2) {\bf 145} (1994), 479-497.
\bibitem{GAR} Garside, F.A., {\it The braid group and other groups}, Quart. J. Math. Oxford Ser. (2) {\bf 78} (1969), 235-254.
\bibitem{WORD1} Garber, D.,Kaplan, S. and Teicher, M., {\it A New Approach for Solving the Word Problem in Braid Groups}, submitted.
\bibitem{EFF} Jacquemard, A., {\it About the effective classification of conjugacy classes of braids}, J. Pure. Appl. Alg. {\bf 63} (1990), 161-169.
\bibitem{BAND} Kang, E.S., Ko, K.H. and Lee, S.J., {\it Band-generator presentation for the 4-braid group}, Top. Appl. {\bf 78} (1997), 39-60.
\bibitem{BGTI} Moishezon, B. and Teicher, M., {\it Braid group techniques in complex geometry I, Line arrangements in $\C \PP^2$}, Contemporary Math. {\bf 78} (1988), 425-555.
\end{\bib}

\end{document}